\documentclass[]{interact}

\usepackage{epstopdf}
\usepackage[caption=false]{subfig}

\usepackage[numbers,sort&compress]{natbib}
\bibpunct[, ]{[}{]}{,}{n}{,}{,}

\theoremstyle{plain}
\newtheorem{theorem}{Theorem}[section]

\newtheorem{corollary}[theorem]{Corollary}

\theoremstyle{definition}

\theoremstyle{remark}

\newcommand{\Cset}{{\mathbb{C}}}

\begin{document}

\articletype{ARTICLE TEMPLATE}

\title{Efficient approximation of functions of some large matrices by partial fraction expansions}

\author{
\name{D. Bertaccini\textsuperscript{a}\thanks{Corresponding author D. Bertaccini. Email: bertaccini@mat.uniroma2.it} \and M. Popolizio\textsuperscript{b} \and F. Durastante\textsuperscript{c}}
\affil{\textsuperscript{a}Universit\`{a} di Roma ``Tor Vergata'', dipartimento di Matematica, viale della Ricerca Scientifica 1, Roma, Italy.
E-mail: bertaccini@mat.uniroma2.it; Istituto per le Applicazioni del Calcolo (IAC) \lq\lq M. Picone\rq\rq, National Research Council (CNR), Roma, Italy. ORCID 0000-0002-3662-278X;\\\textsuperscript{b}Universit\`a del Salento, dipartimento di Matematica e Fisica, via per Arnesano, Monteroni Di Lecce, Lecce, Italy. E-mail: marina.popolizio@unisalento.it; ORCID 0000-0003-0474-2573;\\\textsuperscript{c}Universit\`{a} di Pisa, dipartimento di Informatica, largo Bruno Pontecorvo 3, Pisa, Italy. E-mail: fabio.durastante@di.unipi.it; ORCID 0000-0002-1412-8289.}
}

\maketitle

\begin{abstract}
Some important applicative problems require the evaluation of
functions $\Psi$ of large and sparse and/or \emph{localized}
matrices $A$. Popular and interesting techniques for computing
$\Psi(A)$ and $\Psi(A)\mathbf{v}$, where $\mathbf{v}$ is a vector,
are based on partial fraction expansions. However, some of these
techniques require solving several linear systems whose matrices
differ from $A$ by a complex multiple of the identity matrix $I$ for
computing $\Psi(A)\mathbf{v}$ or require inverting sequences of
matrices with the same characteristics for computing $\Psi(A)$. Here
we study the use and the convergence of a recent technique for
generating sequences of incomplete factorizations of matrices in
order to face with both these issues. The solution of the sequences
of linear systems and approximate matrix inversions above can be
computed efficiently provided that $A^{-1}$ shows certain decay
properties. These strategies have good parallel potentialities.

Our claims are confirmed by numerical tests.
\end{abstract}

\begin{keywords}
matrix functions; partial fraction expansions; large linear systems; incomplete factorizations
\end{keywords}

\begin{amscode}
	65F60, 65F08, 15A23
\end{amscode}

\section{Introduction}
\label{sec:intro}

The numerical evaluation of a function $\Psi(A)\in\Cset^{n\times n}$
of a matrix $A\in\Cset^{n\times n}$ is ubiquitous in models for
applied sciences. Functions of matrices are involved in the solution
of ordinary, partial and fractional differential equations, systems
of coupled differential equations, hybrid differential-algebraic
problems, equilibrium problems, complex networks, in quantum theory,
in statistical mechanics, queuing networks, and many others.
Motivated by the variety of applications, important advances in the
development of numerical algorithms for matrix function evaluations
have been presented over the years and a rich literature is devoted
to this subject; see, e.g.,
\cite{Saad92,Moler.VanLoan.03,Higham2008,HaleHighamTref08} and
references therein.

In this paper we focus mainly on functions of \emph{large} and sparse and/or
\emph{localized} matrices $A$. A typical example of localized matrix
generated by a PDE, say, is the one whose nonnegligible
entries are concentrated in a small region within the computational
domain showing a rapid decay away from this region.
Localization often offers
a way to perform (even full) matrix computations much more
efficiently, possibly with a linear cost with respect to the degrees
of freedom. For a very interesting treatment on this new point of
view we suggest the review~\cite{benzi2016localization}. In the
latter there are also several examples of localized matrices from
physics, Markov chains, electronic structure computations, graph,
network analysis, quantum information theory and many others.

For the computation of $\Psi(A)$ with $A$ as above,
the available literature offers few efficient strategies. The
existing numerical methods for computing matrix functions can be
broadly divided into three classes: those employing approximations
of $\Psi$, those based on similarity transformations of $A$ and
matrix iterations. When the size $n$ of the matrix argument $A$ is
very large, as for example when it stems from a fine grid
discretization of a differential operator, similarity
transformations and matrix iterations can sometimes be not feasible since
their computational cost can be of the order of $n^3$ flops in
general. To overcome these difficulties we consider an efficient
computational framework for approximation algorithms based on
partial fraction expansions. In particular, let us consider an
approximation of $\Psi(A)$ of the form
\begin{equation}\label{eq:sviluppoA}
f(A) =  \sum_{j=1}^{N} c_j (\xi_j I-A)^{-1}
\end{equation}
where scalars $c_j$ and $\xi_j$ can be complex and $I$ is the
$n\times n$ identity matrix.
The above approach has been proven to be effective for a wide set
of functions $\Psi$.

In general, computing \eqref{eq:sviluppoA} requires inverting
several complex valued matrices and, with the exception of lucky or
trivial cases, if $n$ is large, this can be computationally
expensive. We propose to overcome this issue by approximating
directly each term $(\xi_j I-A)^{-1}$ with an efficient update of an
inexact sparse factorization inspired by the complex valued
preconditioners update proposed in~\cite{Bertaccini.04} that there
was defined for symmetric matrices $A$ only. Moreover, such strategy
can be extended to the computation of the action of the matrix
function on vectors, that is, to compute $\Psi(A)\mathbf{v}$ for a
given vector~$\mathbf{v}$. Vectors of this form often represent the
solution of important problems. The simplest example is the vector
$\exp(t_1 A)\mathbf{y}_0$ which represents the solution at a time
$t_1$ of the differential equation $\mathbf{y}^{\prime}(t)=A
\mathbf{y}(t)$ subject to the initial condition
$\mathbf{y}(t_0)=\mathbf{y}_0$.

Note that if the interest is just on obtaining the vector
$\Psi(A)\mathbf{v}$ and not $\Psi(A)$, then \emph{ad hoc} strategies
can be applied as, for example, well known Krylov subspace methods
\cite{Saad92,Hochbruck.Lubich.97} and
\cite{Moret.Novati.04b,vandenEshof.Hochbruck.06,Popolizio.Simoncini.08,AfanasjewEtAl08,KnizhnermanSimoncini2010,MoretPopolizio12,Moret2009,GarrappaPopolizio2011_MCS}
and others.

The paper is organized as follows: in Section
\ref{sec:MatrixFunction} we recall the basics of matrix functions,
together with some results on approximation theory to ground the
proposed approach. Section \ref{sec:update} recalls a recent
updating strategy we propose to use in the algorithms to approximate
matrix functions. In Section \ref{sec:analysis} the proposed
approximation for matrix functions is analyzed by first recalling
some recent results on our updating process for approximate inverse
factorizations and then using the underlying results to build an
a-priori bound for the error made.
Section \ref{sec:Numerical tests} is devoted to numerical tests
showing the effectiveness of the approach in a variety of
applications and comparisons. Section \ref{sec:Discussion} discusses
briefly some final issues.

\section{Computing function of matrices by partial fraction expansions}\label{sec:MatrixFunction}

Many different definitions have been proposed over the years for
matrix functions. We refer to N. Higham~\cite{Higham2008} for an
introduction and references.

In this work we make use of a definition based on the \emph{Cauchy
integral}: given a closed contour $\Gamma$ lying in the region of
analiticity of $\Psi$ and enclosing the spectrum of $A$, $\Psi(A)$
is defined as
\begin{equation}\label{eq:ContourInt}
\Psi(A) =\frac{1}{2\pi i} \int_{\Gamma}  \Psi(z) (zI - A)^{-1} dz.
\end{equation}
Thus, any analytic function $\Psi$ admits an approximation of the
form \eqref{eq:sviluppoA}. Indeed, the application of any quadrature
rule with $N$ points on the contour $\Gamma$, leads to an approximation as in
\eqref{eq:sviluppoA}.

In~\cite{HaleHighamTref08} authors address the choice of the
conformal maps to deal with the contour $\Gamma$ for special
functions like $A^{\alpha}$ and $\log(A)$ when $A$ is a real
symmetric matrix whose eigenvalues lie in an interval $[a,b]\subset
(0,\infty)$. The basic idea therein is to approximate the integral
in \eqref{eq:ContourInt} by means of the trapezoidal rule applied to
a circle in the right half--plane  surrounding $[a,b]$. Thus,
\begin{equation}\label{eq:apprLog}
\Psi(A)\approx f(A)=\gamma A \operatorname{Im}\displaystyle\sum_{j=1}^{N} c_j(\xi_j
I-A)^{-1}
\end{equation}
where $\gamma$ depends on $a, b$ and a complete elliptic integral,
while the $\xi_j$ and $c_j$ involve Jacobi elliptic functions
evaluated in $N$ equally spaced quadrature nodes. We refer to
\cite{HaleHighamTref08} for the implementation details and we make
use of their results for our numerical tests. In particular, an
error analysis is presented there and we report here briefly only
the main result; see \cite{HaleHighamTref08}.
\begin{theorem}\label{thm:3Nquadraturethm}
Let $A$ be a real matrix with eigenvalues in $[a,b]\,,0<a<b$,  let $\Psi$ be
a function analytic in $\Cset\backslash(-\infty,0]$ and let $f(A)$ be the
approximation in \eqref{eq:apprLog}. Then 
$$
\|\Psi(A)-f(A)\| ={\cal{O}}(e^{-{\pi}^2 N/(\log(b/a)+3)}).
$$
\end{theorem}
The analysis in~\cite{HaleHighamTref08} also applies to matrices
with complex eigenvalues.

An approximation like \eqref{eq:sviluppoA} can also derive from a
rational approximation $R_N$ to $\Psi$, given by the ratio of two
polynomials of degree $N$, with the denominator having simple poles.
A popular example is the Chebyshev rational approximation for the
exponential function on the real line. This has been largely used
over the years and it is still a widely used approach, since it
guarantees an accurate result even for low degree $N$, say $N=16$.
Its poles and residues are listed in~\cite{CoMV69} while in
\cite{CaRV84} the approximation error is analyzed and the following
useful estimate is given
$$
\sup_{x\geq 0} |\exp(-x)-R_N(x)|\approx 10^{-N}.
$$
Another example is the diagonal Pad\'e approximation to the logarithm, namely
\begin{equation}\label{eq:Pade}
\log(I+A) \approx f(A) =\displaystyle A\sum_{j=1}^{N} \alpha_j  ( I + \beta_j A)^{-1};
\end{equation}
this is the core of the {\tt logm$\_$pade$\_$pf} code in the package by Higham
\cite{Higham2008} and we will use it in our numerical tests in Section
\ref{sec:Numerical tests}. Unfortunately, as for every Pad\'e approximant,
formula \eqref{eq:Pade} works accurately only when $\|A\|$ is relatively small,
otherwise scaling-and-squaring techniques or similar need to be applied.
The error analysis for the matrix case reduces to the scalar one, according to the
following result; see \cite{KenneyLaub89}.
\begin{theorem}
If $\|A\|<1$ and $f(A)$ is defined as \eqref{eq:Pade} then
$$
\|\log(I+A)-f(A)\|\leq |f(-\|A\|)-\log(1-\|A\|)|.
$$
\end{theorem}

In some important application, the approximation of the matrix
$\Psi(A)$
is not required and it is enough to get the vector $\Psi(A)\mathbf{v}$ 
for a given vector $\mathbf{v}$.
In this case, by using \eqref{eq:sviluppoA}, we formally get the
approximation
\begin{equation}\label{eq:PFE0}
f(A) \mathbf{v} = \sum_{j=1}^{N}c_j(\xi_j I-A)^{-1} \mathbf{v}
\end{equation}
which requires to evaluate $(\xi_j I-A)^{-1}$ or $(\xi_j I-A)^{-1}
\mathbf{v}$ for several values of $\xi_j$, $j=1,\ldots,N$. Usually,
if $A$ is large and sparse or \emph{localized} or even structured,
the matrix inversions in \eqref{eq:PFE0} should be avoided since
each term $\mathbf{w}_{j}\equiv (\xi_j I-A)^{-1} \mathbf{v}$ is
mathematically (but fortunately not computationally) equivalent to
the solution of the algebraic linear system
\begin{equation}
(\xi_j I-A) \ \mathbf{w}_{j} =\  \mathbf{v}.
\label{eq:shifted-systems}
\end{equation}

\section{Updating the approximate inverse factorizations}\label{sec:update}
In the underlying case of interest, i.e., $A$ large and sparse
and/or localized or structured, solving \eqref{eq:shifted-systems}
by standard direct algorithms can be unfeasible and in general a
preconditioned iterative framework is preferable. However, even
using an iterative solver but computing $N$ preconditioners, one for
each of the matrices $(A+\xi_j I)$, can be expensive. At the same
time, keeping the same preconditioner for all the $N$ linear systems
(see, e.g.,~\cite{Popolizio.Simoncini.08}), even if chosen
appropriately, may not account for all the possible issues. Indeed,
very different order of magnitude of the complex valued parameters
$\xi_j$ can cause potential risks for divergence of the iterative
linear system solver. Our proposal is based on cheap updates for
incomplete factorizations developed during the last decade started
by the papers~\cite{Benzi.Bertaccini.03} and~\cite{Bertaccini.04}
essentially based on the inversion and sparsification of a reference
approximation used to build updates. We stress that the updates in
\cite{Benzi.Bertaccini.03} and~\cite{Bertaccini.04} were studied for
symmetric matrices. In recent years these algorithms have been
generalized towards either updates from any symmetric matrix to any
other symmetric (see~\cite{Bertaccini.Sgallari10}) and nonsymmetric
matrices (see
\cite{Bellavia.Bertaccini.Morini10,Bellavia.Bertaccini.Morini11},
and~\cite{BD-interpolated16}) with applications to very different
contexts, but still little attention has been spent on the update of
incomplete factorizations for sequences of nonsymmetric linear
systems with a complex shift.

Among the strategies that can provide a factorization for the
inverse of $A$ we consider the \emph{approximate inverses} or
\emph{AINV} by Benzi et al. (see~\cite{Benzi.survey02} and
references therein) and the \emph{inversion and sparsification}
proposed by van Duin~\cite{vanDuin}, or \emph{INVT} for short. Both
the approaches are very interesting, and differ slightly in their
computational cost (see~\cite{Bertaccini.Filippone12} for some
recent results), parallel potentialities and stability. 

Several efforts have been done in the last decade in order to update
the above mentioned incomplete factorizations in inverse form,
usually as preconditioners; see~\cite{Benzi.Bertaccini.03,Bertaccini.04,Bertaccini.Sgallari10,Bellavia.Bertaccini.Morini11,BD-interpolated16}.

Here, in order to build up an approximate factorization (or, better
saying, to approximate an incomplete factorization)
for each factor $(\xi I-A)^{-1}$, as $\xi$ varies, we assume that
$A$ can be formally decomposed as $A=L\,\,D\,U^H$ with $L$, $U$
lower triangular matrices and that the factorization is well
defined. Then, the inverse of $-A$ can be formally decomposed as
$$
-A^{-1}=U^{-H}D^{-1}L^{-1}= Z D^{-1}  W^H,
$$
where $ W=L^{-H}$ and $Z=U^{-H}$ are upper triangular with all ones
on the main diagonal and $D$ is a diagonal matrix, respectively. The
process can be based also on different decompositions but here we
focus on LDU-types only. In general, this is a not practical way to
proceed because the factors $L$ and $U$ (and thus their inverses)
are often dense. At this point we have two possibilities. The first
is use AINV and its variants (again see \cite{Benzi.survey02}) that
provides directly an approximate inverse in factored form for $A$
whose factors can be suitably sparse as well if $A$ is sparse or
shows certain decay properties. The second is use an inversion and
sparsification process as in \cite{vanDuin}, that, starting from a
sparse incomplete factorization for $A$ such as ILU (see, e.g.,
\cite{saad2003iterative}) $P=\tilde L \tilde D \tilde U^H$
approximating $A$, whose factors $\tilde L$, $\tilde U$ are sparse,
produces an efficient inversion of $\tilde L$, $\tilde U$ and
provides also a post-sparsification of the factors $Z$ and $W$ to
get $\tilde Z$ and $\tilde W$. A popular post-sparsification
strategy can be to zero all the entries smaller than a given value
and/or outside a prescribed pattern. We call \emph{seed
preconditioner}, denoted $P_0$, the following approximate
decomposition of $A^{-1}$:
\begin{equation}\label{eq:PinvSeed}
P_0=\tilde Z\, \tilde D^{-1}\, \tilde W^H.
\end{equation}
Similarly to what done above for $A^{-1}$, in the style
of~\cite{Bertaccini.04}, given a complex pole $\xi$, a factorization
for the inverse of the complex nonsymmetric matrices in
\eqref{eq:shifted-systems} can be formally obtained by the
identities
\begin{eqnarray}\label{eq:Axi}
A_{\xi}^{-1}&\equiv& (-A+\xi\,I)^{-1}  \\ \nonumber &=&
(W^{-H}\,D\,Z^{-1}+\xi \, W^{-H} (W^H\,Z) Z^{-1})^{-1} \\
\nonumber &=& Z \, \left( D+\xi\, E \right)^{-1} W^H, \quad
E=W^H\,Z.
\end{eqnarray}
However, as recalled above, the factors $Z$ and $W$ are dense in general. Therefore,
their computation and storage are sometimes possible for $n$ small to moderate
but can be too expensive to be feasible for $n$ large. This issue can be faced by
using the sparse approximations $\tilde Z$ and $\tilde W$ for $Z$ and $W$, respectively,
produced by AINV, by inversion and sparsification or by another process generating a sparse factorization for $A^{-1}$. Indeed, supposing that the chosen algorithm generates
a well defined factorization,  we can provide an approximate
factorization for the inverse of $A+\xi_j I$. 
In particular, we get a sequence of approximate factorization
candidates using $P_0$ defined above as a reference and $\tilde E$, a
sparsification of the nonsymmetric \emph{real valued} matrix $E$,
with the approximation of $A_{\xi}^{-1}$ given by $P_{\xi}$ defined
as
\begin{equation}\label{eq:Pk}
P_{\xi}=\tilde Z \, \left(\tilde D +\xi \tilde E\right)^{-1}\tilde
W^H,
\end{equation}
where, by using the formalism introduced in
\cite{Bellavia.Bertaccini.Morini11},
\begin{equation}\label{eq:banded}
\tilde E=g(\tilde W^H\,\,\tilde Z).
\end{equation}
The function $g$ serves to generate a sparse matrix from a full one
such that the linear systems with matrix $\tilde D +\xi \tilde E$
can be solved with a low computational complexity, e.g., possibly
linear in $n$.
As an example, if the entries of $A^{-1}$ decay fast away from
the main diagonal,
we can consider the sparsifying function $g=g_m$, 
$$
g_m:\mathbb{C}^{n\times n}\rightarrow \mathbb{C}^{n\times n},
$$
extracting $m$ upper and lower bands (with respect to the main
diagonal, which is the $0$-diagonal) of its matrix argument
generating an $(2m,2m)$--\emph{banded matrix}. 
In general, a matrix $A$ is called
$m$-\emph{banded} if there is an index $l$ such that
\begin{equation*}
a_{i,j} = 0, \;\;\text{ if }j \notin [i-l,i-l+m].
\end{equation*}
It is said to be \emph{centered} and $m$-\emph{banded} if $m$ is
even and the $l$ above can be chosen to be $m/2$. In this case the
zero elements of the \emph{centered} and $(m,m)$-\emph{banded} are:
\begin{equation*}
a_{i,j} = 0, \;\;\text{ if }|i-j|>\frac{m}{2},
\end{equation*}
thus selfadjoint matrices are naturally centered, i.e., a tridiagonal selfadjoint matrix is centered and $2$-banded.
This choice will be used
in our numerical examples but of course different choices for $g$ can
be more appropriate in different contexts.
A substantial saving
can be made by approximating
$$
g_m(\tilde W^H\,\,\tilde Z) \textrm{ with } g_m(\tilde
W^H)\,g_m(\,\tilde Z), \quad m>0
$$
with a reasonable quality of the approximation, i.e., under suitable
conditions and provided $m>0$, the relative error
$$
\frac{||g_m(\tilde W^H\,\,\tilde Z)-g_m(\tilde W^H)\,g_m(\,\tilde Z)
||}{||\tilde W^H\,\,\tilde Z||}
$$
can be moderate in a way that will be detailed in Theorem~\ref{teo:gm(Z)=g(Z)g(W)} discussed in the next section.

\section{Analysis of the approximation processes}
\label{sec:analysis}

We use here the underlying approximate inverses in factored form
\eqref{eq:Pk} as a preconditioner for Krylov solvers to approximate
$\Psi(A) \mathbf{v}$ and to approximate $f(A)$ by
\begin{align}\label{eq:f_tilde}
\tilde f(A) = & \displaystyle\sum_{j=1}^{N} c_j P_{\xi_j} \\
= & \displaystyle\sum_{j=1}^{N} c_j \tilde Z \, \left(\tilde D +\xi_j
\tilde E\right)^{-1}\tilde W^H.\nonumber
\end{align}
In order to discuss an a-priori bound for the norm of the error
$||\Psi(A)-\tilde f(A)||$ generated by the various approximation
processes, supposing we are operating in exact arithmetic, we need
some results on the update of the approximate inverse
factorizations.


Let us recall a couple of results that can be derived as corollaries
of Theorem~4.1 in~\cite{Demko.Moss.Smith84}. 
In this context, we consider a general complex, separable, Hilbert space $H$, and denote
with $\mathcal{B}(H)$ the Banach algebra of all linear operators on
$H$ that are also bounded. If $A \in \mathcal{B}(H)$, then $A$ can
be represented by matrix with respect to any complete orthonormal
set thus $A$ can be regarded as an element of $B(l^2(S))$, a matrix
representing a bounded operator in $\mathcal{B}(l^2(S))$, where $S =
\{1,2,\ldots,N\}$. 
\begin{theorem}
\label{teo:decay} Let $A\in \mathbb{C}^{n\times n}$ be a nonsingular
$(2m,2m)$--banded matrix, with $A \in \mathcal{B}(l^2(S))$ and
condition number $\kappa_2(A)\geq 2$. Then, by denoting with
$b_{i,j}$ the ${i,j}-$entry of $A^{-1}$ and with
\begin{equation*}
\beta=\left( \frac{\kappa_2(A)-1}{\kappa_2(A)+1} \right)^{\frac{1}{2m}},
\end{equation*}
for all $\tilde\beta>\beta$, $\tilde \beta<1$, there exists a constant $c>0$ such that
\begin{equation*}
|b_{i,j}|\leq c\, \tilde \beta^{|i-j|},
\end{equation*}
with
\begin{equation*}
\begin{split}
c\leq & (2m+1)\frac{\kappa_2(A)+1}{\kappa_2(A)-1} \|A^{-1}\|\kappa_2(A) \\
\leq & 3(2m+1)\|A^{-1}\|\kappa_2(A).
\end{split}
\end{equation*}
\end{theorem}
For the proof and more details, see \cite[Theorem 3.10]{DanieleFabioBook2018}.


We can note immediately that the results in Theorem \ref{teo:decay},
without suitable further assumptions, can be of very limited use
because:
\begin{itemize}
\item the decay of the extradiagonal entries can be very slow, in
principle arbitrarily slow;
\item the constant $c$ in front of the bound depends on the
\emph{condition number} of $A$ and we are usually interested in
approximations of $A^{-1}$ such that their condition numbers can range from
moderate to high;
\item the bound is far to be tight in general. A trivial example
is given by a diagonal matrix with entries $a_{j,j}=j$,
$j=1,\ldots,n$. We have that $b_{i,j}=0$, $i\neq j$ but of course
$\kappa_2(A)=a_{n,n}/a_{1,1}=n$.
\item If we take $m=n$ and $n$ is very large, then $\tilde
\beta$ must be chosen very near $1$ and it is very likely
that no decay can be perceptible with the bound in Theorem \ref{teo:decay}.
\end{itemize}
However, the issues presented here are more properly connected with
the decay properties of the matrices $Z$, $W$ (and therefore $\tilde
Z$, $\tilde W$). Using similar arguments as in Theorem 4.1 in
\cite{Benzi.Tuma.decay00}, it is possible to state the following
result.
\begin{corollary}
\label{coro:decayZW}
Let $A\in \mathbb{C}^{n\times n}$ be invertible, $A \in \mathcal{B}(l^2(S))$,
and with its symmetric part positive definite.
Then for
all $i$, $j$ with $j > i$, the entries $z_{i,j}$ in $Z = L^{-H}$ and
$w_{i,j}$ in $W = U^{-1}$ satisfy the following upper bound:
\begin{equation*}
|z_{i,j}|\leq c_1\, \tilde \beta_1^{j-i}, \quad |w_{i,j}|\leq c_2\,
\tilde \beta_2^{j-i}, \quad j>i
\end{equation*}
(note that $z_{i,j}$, $w_{i,j}=0$ for $j\leq i$), where
\begin{equation*}
0< \tilde \beta_1, \  \tilde \beta_2 \leq \tilde \beta < 1
\end{equation*}
and $c_1$, $c_2$ are positive constants, $c_1,c_2\leq c_3\cdot
\kappa_2(A).$
\end{corollary}

Recently, this kind of decay bound for the inverses of matrices was
intensely studied, and appears also with other structures. Consider,
e.g., the case of nonsymmetric band matrices in
\cite{nabben1999decay}, tridiagonal and block tridiagonal matrices
in~\cite{meurant1992review}, triangular Toeplitz matrices coming
from the discretization of integral equations~\cite{ford2014decay},
Kronecker sum of banded matrices~\cite{canuto2014decay}, algebras
with structured decay~\cite{jaffard1990proprietes} and many others.
Thus, the results we propose can be readily extended to the above
mentioned cases.

If the seed matrix $A$ is, e.g., diagonally dominant, then the decay
of the entries of $A^{-1}$ and therefore of $W$, $Z$ ($\tilde W$,
$\tilde Z$) is faster and more evident. This can be very useful for
at least two aspects:
\begin{itemize}
\item the factors  $\tilde W$,
$\tilde Z$ of the underlying approximate inverse in factored form
can show a narrow band for drop tolerances even just slightly
larger than zero;
\item banded approximations can be used not only for post--sparsifying
$\tilde W$, $\tilde Z$ in order to get more sparse factors, but
also the update process can benefit from the fast decay.
\end{itemize}

\begin{theorem}
\label{teo:gm(Z)=g(Z)g(W)} Let $A\in \mathbb{C}^{n\times n}$ be
invertible, $A \in \mathcal{B}(l^2(S))$, and with its symmetric part positive definite.
Let $g_m = [\cdot]_m$ be a sparsifying function extracting the $m$
upper and lower bands of its argument. Then, given the matrices from Corollary~\ref{coro:decayZW}, we have
$$
[\tilde W^H \, \tilde Z]_m=[\tilde W^H]_m[\tilde Z]_m+R(A,m),
$$
$$
|(R(A,m))_{i,j}|\leq c_4\tilde \beta^{|i-j|},
$$
where $ c_4 =  c_1 c_2$.
\end{theorem}
The above result can be proved by comparing the expressions of
$\tilde W^H \tilde Z$ and $g_m(\tilde W^H \tilde Z)$ and using
Corollary \ref{coro:decayZW} with an induction argument on $n$.

As a matter of fact, we see that a fast decay of entries of $A^{-1}$
guarantees that the essential component of the proposed update
matrix, i.e.,  $\tilde E=\tilde W^H  \tilde Z$, can be cheaply,
easily and accurately approximated by the product $g_m(\tilde
W^H)\cdot g_m(\tilde Z)$, without performing the time and memory
consuming matrix-matrix product $\tilde W^H \tilde Z$.

On the other hand, if the decay of the entries of $A^{-1}$ is fast,
even a simple diagonal approximation of $\tilde W^H \, \tilde Z$ can
be accurate enough. In this case, there is no need to apply the
approximation in Theorem \ref{teo:gm(Z)=g(Z)g(W)}. The update matrix
$\tilde E$ can be produced explicitly by the exact expression of
$\operatorname{diag}(\tilde W^H \, \tilde Z)$ we give in the following corollary.
%
\begin{corollary}
\label{coro:g1(Z)diretto} Let $A\in \mathbb{C}^{n\times n}$ be
invertible, $A \in \mathcal{B}(l^2(S))$, and with its symmetric part
$(2m,2m)$--banded and positive definite, $1\leq m\leq n$. Then, the
diagonal approximation for $\tilde E$ generated by the main diagonal
of $\tilde W^H \, \tilde Z$ is given by
\begin{equation*}
\tilde E=\operatorname{diag}(\tilde W^H \, \tilde Z)=(d_{i,i}),
\end{equation*}
where
\begin{equation*}
d_{i,i}=1+\sum_{j=1, \ i-j\leq m}^{i-1} \tilde w_{j,i} \tilde
z_{j,i}, \quad 1\leq i\leq n,
\end{equation*}
where $\tilde w_{j,i}$ and $z_{j,i}$ are the entries of $\tilde W$
and $\tilde Z$, respectively.
\end{corollary}

All our numerical experiments use the approximations proposed in
Theorem~\ref{teo:gm(Z)=g(Z)g(W)} and in Corollary
\ref{coro:g1(Z)diretto} without perceptible loss of accuracy; see
Section \ref{sec:Numerical tests}.

We use the properties stated by the previous results to get an
\emph{a--priori} estimates of the global error
$\|\psi(A)-\tilde{f}(A)\|$ that is given below in exact arithmetics
and for $A$ symmetric and definite positive in order to use the
results in \cite{HaleHighamTref08}.
Note that, with the above hypotheses, we get $\tilde Z=\tilde W$
in~\eqref{eq:f_tilde} and therefore
\begin{equation}\label{eq:f_tilde_Asymm}
\tilde f(A) = \sum_{j=1}^{N} c_j \tilde Z \, \left(\tilde D +\xi_j
\tilde E\right)^{-1}\tilde Z^H.
\end{equation}

\begin{theorem}
\label{teo:error-bound} Let $A$ be a real positive definite matrix
with eigenvalues in $[a,b]$, $0<a<b$, $\Psi$ a function analytic in
$\Cset\backslash(-\infty,0]$, $f(A)$ the approximation of $\psi(A)$
in~\eqref{eq:apprLog} and $\tilde{f}(A)$ the approximation of $f(A)$
in (\ref{eq:Pk}) with $0<\tau<1$ drop tolerance used to produce
$\tilde Z$. Moreover, let $\Delta_Z=\tilde Z-Z$. Then,
\begin{equation*}
\|\psi(A)-\tilde{f}(A)\| \leq E_1(N) + E_2(\tau),
\end{equation*}
with
\begin{equation*}
E_2(\tau) = c(\beta)\,\tau=O(||\Delta_Z||),
\end{equation*}
where $c=c(\beta)$ is a parameter that depends on the decay of the offdiagonal entries of
$Z$, and
\begin{equation*}
E_1(N) ={\cal{O}}(e^{-{\pi}^2 N/(\log (b/a)+3)}).
\end{equation*}
\end{theorem}
The above result follows by observing that
\begin{eqnarray}
\|\psi(A)-\tilde{f}(A)\| =  & \|\psi(A)-f(A) + f(A)-\tilde{f}(A)\| \label{eq:psi-tildef}\\
\leq & \|\psi(A)-f(A)\| + \|f(A)-\tilde{f}(A)\|\nonumber.
\end{eqnarray}
The upper bound for the quadrature error $E_1(N) = \|\psi(A)-f(A)\|$
for an analytic function $\Psi$ is obtained straightforwardly from
Theorem~\ref{thm:3Nquadraturethm}. Recall that the bound is derived
from the classical error estimate for the Trapezoidal/Midpoint rule
\cite[Section~4.6.5]{davis2007methods}. Similar bounds for functions
that are less smooth can be provided as well, even if they show just
a polynomial decay, see again \cite[Section~2.9]{davis2007methods}.

The upper bound for the errors generated by the approximation of the
terms $(A+\xi I)^{-1}$ by the approximate inverse factorization
updates, i.e., $E_1(N) = \|f(A)-\tilde{f}(A)\|$, is easily derived
by working on the norm of the difference between (\ref{eq:Axi}) and
(\ref{eq:Pk}) substituting to $\tilde Z$ the expression $\tilde
Z=Z+\Delta_Z$ and to $(\tilde D+\xi \tilde E)^{-1}$ the expression
$(D+\xi E)^{-1}+\Delta_D$. The claim follows by observing that
$||\Delta_D||$, $||\Delta_Z||$ can be bounded by $c(\beta)\,\tau$;
see Theorem  \ref{teo:decay} and Corollary \ref{coro:decayZW}.

A generalization of Theorem \ref{teo:error-bound} for nonsymmetric
matrices $A$ can be given with similar arguments.

The main purpose of the a-priori upper bound in Theorem
\ref{teo:error-bound} should be intended as more qualitative than
quantitative, for showing that the multiple approximation processes
considered here for computing $\tilde{f}(A)$ \emph{converge}, i.e.,
in exact arithmetic, under the hypotheses of Theorem
\ref{thm:3Nquadraturethm}, $\|\psi(A)-\tilde{f}(A)\|\rightarrow 0$
if $\tau\rightarrow0$ and $N\rightarrow\infty$.

%

\subsection{Cross-relations between the function $g$ and drop tolerance $\tau$}
To clarify the role of the function $g$ introduced in
\eqref{eq:banded} and the drop tolerance $\tau$ for AINV, we compare
the results of our approach to compute $\exp(A)$ with the built--in
Matlab function {\tt expm}. We use the expression in
\eqref{eq:sviluppoA} for the Chebyshev rational approximation of
degree $N=16$ so that we can consider the approximation error
negligible.

\begin{figure}[htb]
\centering
\subfloat{\includegraphics[width=0.5\columnwidth]{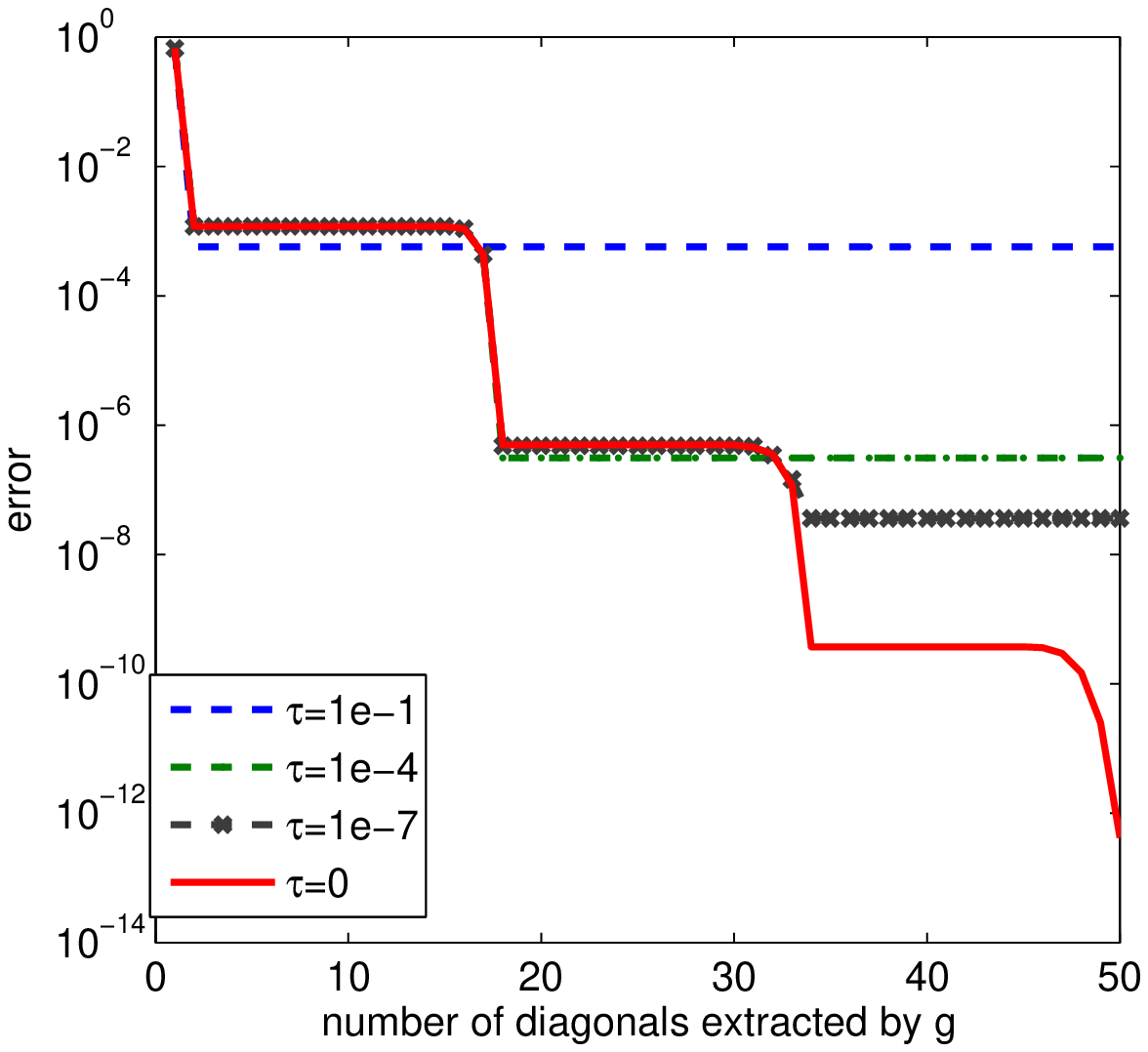}}
\subfloat{\includegraphics[width=0.5\columnwidth]{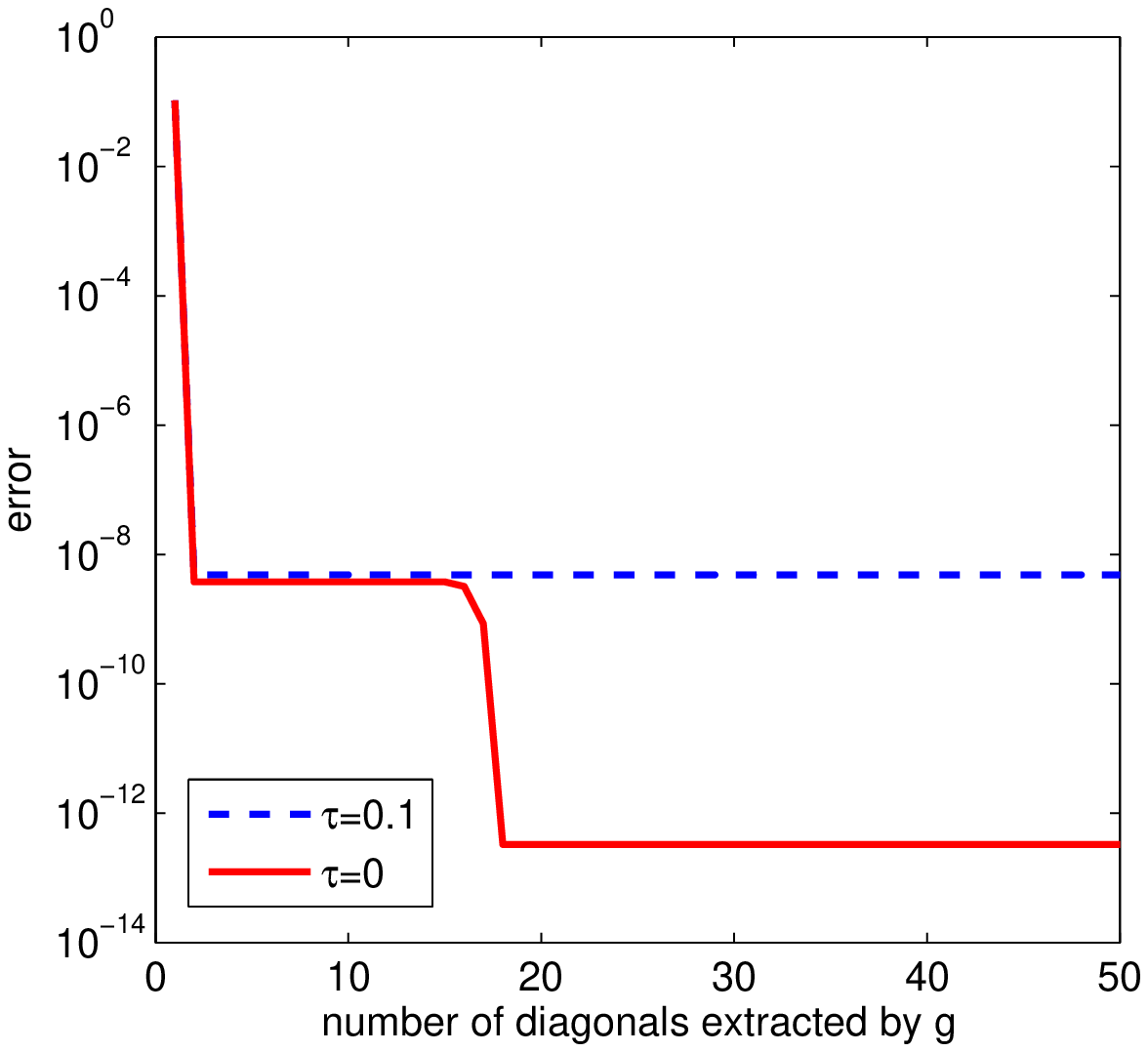}}
\caption{Behavior of the error for $\exp(A)$ as $\tau$ and $g$ vary.
The $50\times 50$  matrix argument $A$ has the expression in
\eqref{eq:TestMat} with $\alpha=\beta=0.5$ (left),
$\alpha=\beta=1.2$ (right). The  x-axis reports the number of
diagonals the function $g$ selects while the y-axis reports the
error with respect to the Matlab's {\tt{expm(A)}}. AINV is used with the tolerance $\tau$ given in the figures' caption.}
\label{ErrorVsBand}
\end{figure}

Consider the localized matrix $A=(A_{i,j})$ described in~\cite{Benzi.Razouk.07}
with entries
\begin{equation}\label{eq:TestMat}
A_{i,j}= \left\lbrace\begin{array}{ll} e^{-\alpha(i-j)}, &  i\geq j,\\
e^{-\beta(j-i)} & i<j,
\end{array}\right.\qquad \alpha,\beta>0.
\end{equation}
This is typical example of a localized matrix, completely dense but
with rapidly decaying entries. These matrices are usually replaced
with banded matrices obtained by considering just few bands or by
dropping entries which are smaller than a certain threshold. Here we
sparsify it by keeping only $15$ off--diagonals on either side of
its main diagonal. Let us take $\alpha=\beta=0.5$ and
$\alpha=\beta=1.2$ for a small example, namely $50\times 50$, in
order to show the application of Theorem~\ref{teo:decay}. The
approximation we refer to is \eqref{eq:f_tilde}, in which we let
$\tau$ and  $g$ change, with effects on the factors $\tilde Z,
\tilde W$ and  $\tilde E$ in \eqref{eq:banded}, respectively. The
continuous curves in Figure \ref{ErrorVsBand} refer to the ``exact''
approach, that is, for $\tau=0$ leading to full factors $\tilde W$
and $\tilde Z$. In the abscissa we report the number of
extra-diagonals selected by $g$. Notice that both $\tau$ and $g$ are
important because even for $\tau=0$ more extra-diagonals are
necessary to reach a high accuracy.

From the plots in Figure \ref{ErrorVsBand}, we note that the loss of
information in discarding entries smaller than $\tau$ cannot be
recovered even if $g$ extracts a full matrix. In the left plot, for
a moderate decay in the off-diagonals entries, a conservative $\tau$
is necessary to keep the most important information. On the other
hand, when the decay is more evident, as in the right plot, a large
$\tau$ is enough, and $g$ keeping just two diagonals gives already a
reasonable accuracy. We get similar results also for the logarithm,
as well as for other input matrices.

\subsection{Choosing the reference preconditioner(s)}
\label{sec:seed-precond}

\begin{figure}[htb]
\centering
\subfloat{\includegraphics[width=0.5\columnwidth]{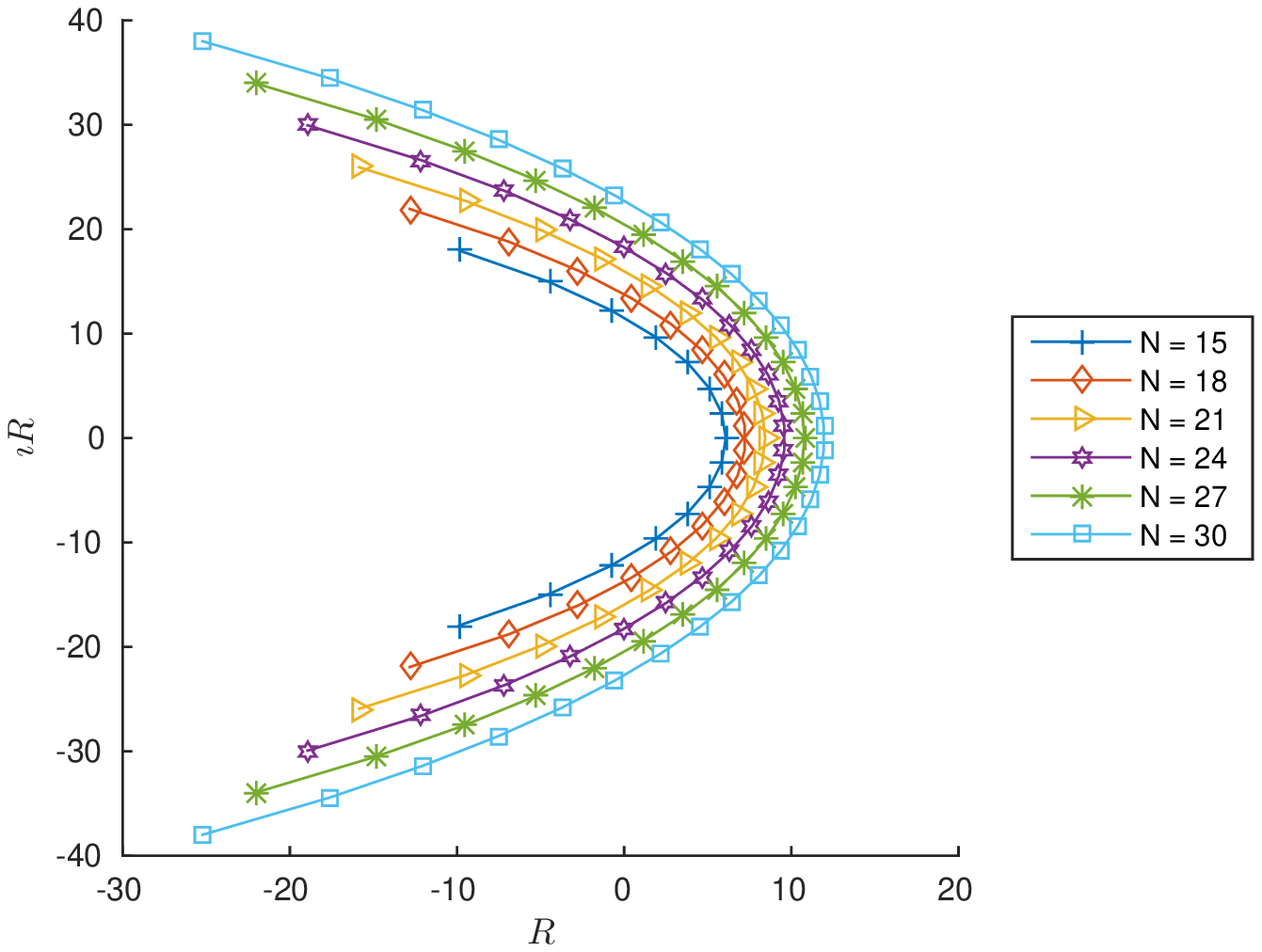}}
\subfloat{\includegraphics[width=0.5\columnwidth]{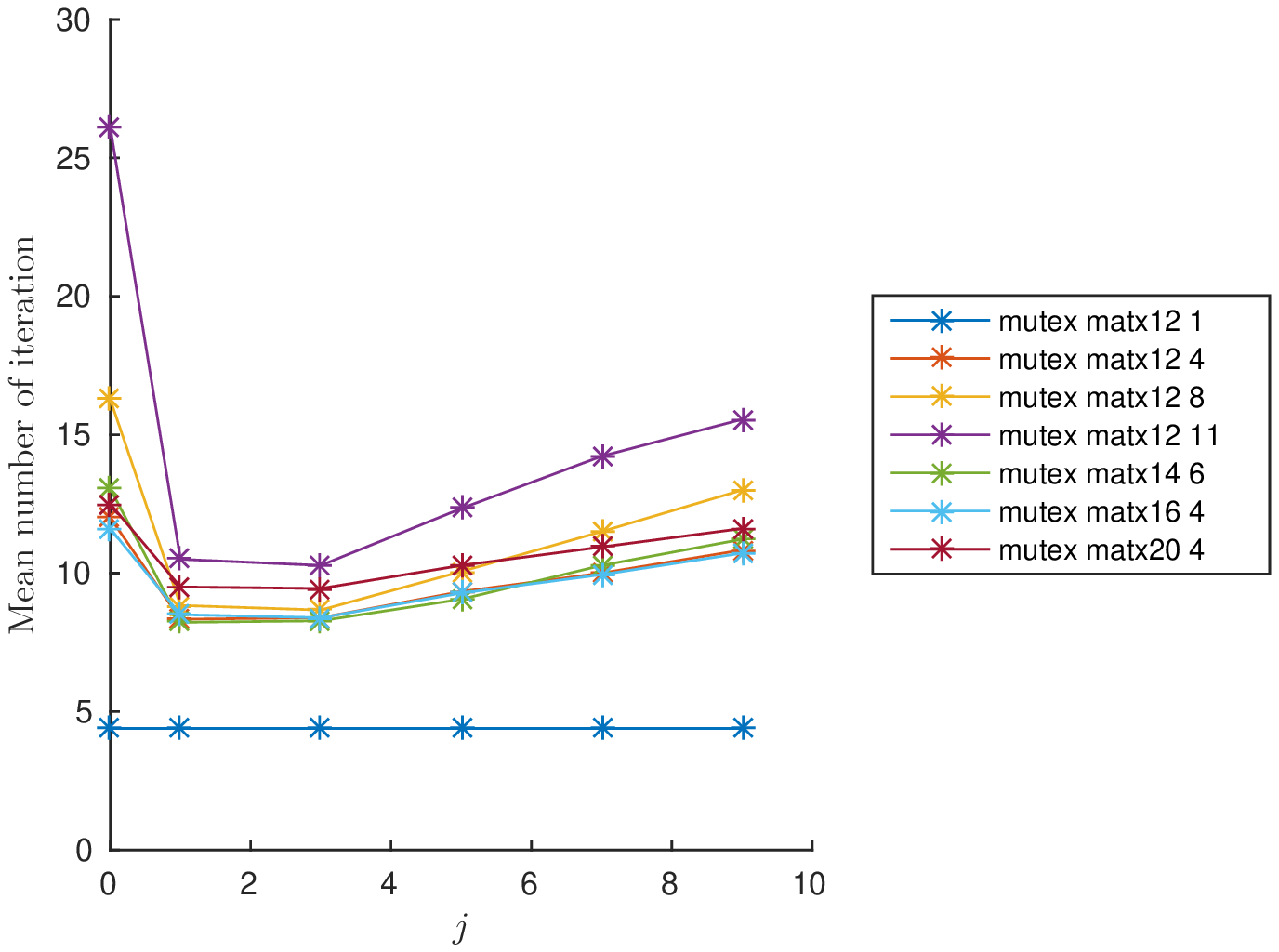}}
\caption{Position of the poles of Chebyshev approximation of
$\exp$ (left) and a sample of mean number of iterations (right) for different choice of
$P_{\text{seed}}$ for the \emph{mutual exclusion model} from~\cite{william1991marca}.}
\label{fig:polesofapprox}
\end{figure}

To generate a viable update \eqref{eq:Pk}, we need to compute an
appropriate seed preconditioner \eqref{eq:PinvSeed}. Note that the
poles $\xi_j$ in the partial fraction expansion \eqref{eq:PFE0} for
the Chebyshev approximation of the exponential have a modulus that
grows with the number of points; see, e.g., Figure~\ref{fig:polesofapprox}.
Therefore, we need to take into account the possibility that the
resolvent matrices
\begin{equation*}
(\xi_j I - A)^{-1}, \quad j=1,\ldots,N
\end{equation*}
become diagonally dominant or very close to the identity matrix, up
to a scalar factor, or, in general, with a spectrum that is far from
the one of $-A$. Sometimes the matrices related to the resolvent
above can be so well conditioned that the iterative solver
does not need any preconditioner. 
In this case, any choice of the seed preconditioner as an
approximate inverse of the matrix $-A$ is almost always a poor
choice, and thus also the quality of the updates; see
\cite{Benzi.Bertaccini.03,Bertaccini.04}. Let us consider, e.g., the
matrices from the \emph{mutual exclusion model}
from~\cite{william1991marca}. These are the transition matrices for
a model of $M$ distinguishable process (or users) that share a
resource, but only $M^{\prime}$, with $1 \leq M^{\prime} \leq M$
that could use it at the same time. In
Figure~\ref{fig:polesofapprox} (right) the underlying experiments
are reported as ``{mutex matx$M$ $M^{\prime}$}''. We report  the
mean iterations required when the $P_{\text{seed}}$ corresponding to
$\xi_j$ is used for $j=1,\ldots, N$, while $j=0$ refers to the seed
preconditioner for $-A$, all obtained with INVT for $\tau_L = 1e-5$
and $\tau_Z = 1e-2$. The plot clearly confirms that working with
$-A$ is always the most expensive choice, while better results are
obtained for whatever pole and sometimes the pole with the largest
modulus slightly betters the others.

Observe also that in this way complex arithmetic should be used to
build the approximate inverse of the matrix $(\xi_1 I - A)$, because
its main diagonal has complex valued entries.

\section{Numerical tests}
\label{sec:Numerical tests}

The codes are written in Matlab (R2016a). The machine used is a laptop running Linux
with 8Gb memory and CPU Intel(R) Core(TM) i7-4710HQ CPU with clock
2.50GHz.

The sparse inversion algorithm chosen for each numerical test (those
used here are described in Section \ref{sec:update}) takes into
account the choice made for the computation of the reference (or
\emph{seed} for short) preconditioners. If the matrix used to
compute the seed preconditioner is real, we use the AINV. Otherwise,
the inversion and sparsification of the ILUT Algorithm, or
\emph{INVT} for short, requiring a dual threshold strategy. See
\cite{Bertaccini.Filippone12} for details and a revisitation of AINV
and INVT techniques. 
In the
following, the symbols
$\tau$ denotes drop tolerance for AINV while $\tau_{L}$, $\tau_{Z} $
the threshold parameter for ILU decomposition and for
post-sparsification of the inverted factors of INVT; respectively;
see also the details discussed in Section~\ref{sec:seed-precond}.
$\varepsilon_{\text{rel}}$ denotes the standard relative (to a
reference solution) error.

Other details on the parameters and strategies used are given in the
description of each experiment.

\subsection{Approximating $\Psi(A)$}

Let us focus on the approximation of $\exp(A)$ and $\log(A)$. In the
following tables, the columns {\emph{Update}} refers to the approximation \eqref{eq:f_tilde}.
Columns {\emph {Direct}} are based on the
direct inversion of the matrices $(A+\xi_j I)^{-1}$ in \eqref{eq:sviluppoA}.

The Fill--In for computing the incomplete factors approximating the
underlying matrices is computed as
\begin{equation}\label{FI}
\text{Fill--In} \ =\frac{nnz(\tilde Z)+nnz(\tilde W) -n}{n^2},
\end{equation}
where $n$ denotes the size and $nnz(\cdot)$ the number of the
nonzero entries, as usual.

\subsubsection{$\log(A)$ - Exponential Decay} We consider the evaluation of $\log(A)$ where the entries of $A$ are
as in \eqref{eq:TestMat} with $\alpha=0.2$ and $\beta=0.5$ and $n$
varies from $500$ to $8000$. For this matrix we use a {drop
tolerance} $\tau=0.1$ to get a sparse approximate inverse
factorization of $A$ with AINV. The resulting factors $\tilde{Z}$ and
$\tilde{W}$ are bidiagonal and thus we take $g(X)=X$. The inversion
of the tridiagonal factors is the more demanding part of the
{\emph{Update}} technique. For this test, we compare the
{\emph{Update}} and {\emph{Direct}} methods, based on the
approximation \eqref{eq:apprLog}, with the Matlab function {\tt
logm} and the {\tt logm$\_$pade$\_$pf} code in the package by N.
Higham~\cite{Higham2008}.

Numerical tests on scalar problems show that the degree $N=5$ for
the Pad\'e approximation \eqref{eq:Pade} and $N=7$ for the
approximant in \eqref{eq:apprLog} allow to reach a similar accuracy
with respect to the reference solution. Thus, we use these values
for $N$ in our tests.

\begin{table}[htb]
\centering
\caption{Execution time in seconds for $\log(A)$ for $A$ as in \eqref{eq:TestMat}
with $\alpha=0.2$ and $\beta=0.5$. AINV with $\tau=1e-1$ is used. 
}\label{Tab:log}
\begin{tabular}{rrrrrr}
\toprule
n  & \textrm{Update} & \textrm{Direct} & \textrm{logm} &
\textrm{logm$\_$pade$\_$pf} & \textrm{Fill--In}\\
\midrule
500&  1.53 &  1.33 & 13.05 &  {\bf{0.67}} &6e-3\\
1000&  4.90 &  4.69 & 44.40 &  {\bf{3.31}}&3e-3\\
2000& {\bf{12.23}} & 13.86 & 407.28 & 38.67&1e-3\\
4000& {\bf{37.04}} & 56.23 & 6720.36 & 522.25&7e-4\\
8000& {\bf{168.41}} & 412.30 & 70244.41 & 6076.00&7e-4\\
\bottomrule
\end{tabular}
\end{table}

Results in Table \ref{Tab:log} show that, for small  examples, the
{\emph{Update}} and the {\tt{logm$\_$pade$\_$pf}} approaches require
a similar execution time, while the efficiency of the former becomes
more striking with respect to all the others as the problem
dimension increases.

\subsubsection{$\exp(A)$ - Exponential Decay} We now consider the error for the matrix exponential. The test
matrix is symmetric as in \eqref{eq:TestMat} for three choices of
the parameter $\alpha$. We analyze the error of the approximations
provided by the {\emph{Update}} and  {\emph{Direct}} methods, for
the Chebychev rational approximation of degree $N=8$, with respect
to the results obtained by the {\tt{expm}} Matlab command. We
consider $\alpha=\beta=1,\,\alpha=\beta=1.5,\,\alpha=\beta=6$. For
the first two cases, the drop tolerance for the AINV is $\tau=0.1$
and $g=g_1(\cdot)$ extracts just the main diagonal and one
superdiagonal. For the third case, AINV with $\tau=10^{-3}$ is used
and $\tilde Z,\,\tilde W$ are both diagonal. No matrix inversion is
thus performed.

\begin{table}[htb]
\centering \normalsize \caption{Errors for the {\emph{Update}} and
{\emph{Direct}} methods compared to the Matlab's {\tt{expm(A)}}. The
parameters are $\tau=0.1$, $\tilde Z$ and $\tilde W$ bidiagonal, for
$\alpha=\beta=1$ (left) and $\alpha=\beta=1.5$ (center);
$\alpha=\beta=6$. AINV is used with $\tau=10^{-3}$ and $\tilde Z$,
$\tilde W$ are diagonal (right).}\label{tab:errorExpMat}
\subfloat{
\begin{tabular}{ccc}
\toprule
n  & \text{Update} & \text{Direct}   \\
\midrule
500&  1.1e-7 &  2.3e-8 \\
1000&  1.1e-7 &  2.3e-8 \\
2000&  1.1e-7 &  2.3e-8 \\
4000&  1.1e-7 &  2.3e-8 \\
\bottomrule
\end{tabular}
}\hfil
\subfloat{
\begin{tabular}{ccc}
\toprule
n  & \text{Update} & \text{Direct} \\
\midrule
500 & 2.3e-08 & 2.3e-08 \\
1000 & 2.3e-08 & 2.3e-08 \\
2000 & 2.3e-08 & 2.3e-08 \\
4000 & 2.3e-08 & 2.3e-08 \\
\bottomrule
\end{tabular}
}\hfil
\subfloat{
\begin{tabular}{ccc}
\toprule
n  & \text{Update} & \text{Direct}  \\
\midrule
500 & 4.5e-06 & 1.8e-08 \\
1000 & 4.5e-06 & 1.8e-08 \\
2000 & 4.5e-06 & 1.8e-08 \\
4000 & 4.5e-06 & 1.8e-08 \\
\bottomrule
\end{tabular}
}
\end{table}

Results from Table \ref{tab:errorExpMat} show the good accuracy the
{\emph{Update}} approach reaches. Indeed, although the presence of
the sparsification errors (see the action of $\tau$ and $g$), the
error is comparable with the one of the {\emph{Direct}} method,
which does not suffer from truncation.
For the case $\alpha=\beta=6$, the difference between the two errors
is more noticeable but it has to be balanced with great savings in
timings. Indeed, in this case the decay of the off--diagonal entries
of the inverse of $A$ is very fast and we exploit this feature by
combining the effect of the small drop tolerance $\tau=10^{-3}$ and
a function $g=g_0(\cdot)$ extracting just the main diagonal. Then,
the computational cost is much smaller for the {\emph{Update}}
approach since no matrix inversion is explicitly performed and we
experienced an overall linear cost in $n$, as in the other
experiments. Thus, when a moderate accuracy is needed, the
{\emph{Update}} approach is preferable, since it is faster; see
Table \ref{tab:timeExpDiag}.

\begin{table}[htb]
\centering \caption{Timings in seconds for $\exp(A)$ with $A$ as in
\eqref{eq:TestMat} with $\alpha=\beta=6$. AINV with $\tau=10^{-3}$
is used and $g$ extracts just the main diagonal. The
{\emph{Fill--In}} column refers to the Fill--In occurred for
computing the factors $\tilde W$ and $\tilde Z$ measured as in
\eqref{FI}.}\label{tab:timeExpDiag}
\begin{tabular}{rrrrr}
\toprule
n  & \textrm{Update} & \textrm{Direct} & \tt{expm} & \textrm{Fill--In}\\
\midrule
500&  \textbf{0.01} &  0.06 &  1.97 & 2.0e-3\\
1000&  \textbf{0.00} &  0.01 &  6.19 & 1.0e-3\\
2000&  \textbf{0.00} &  0.01 & 30.52 & 5.0e-4\\
4000&  \textbf{0.00} &  0.06 & 172.89&2.5e-4\\
8000&  \textbf{0.01} &  0.10 & 910.16 & 1.3e-4\\
\bottomrule
\end{tabular}
\end{table}

%
\subsubsection{$\exp(A)$ - Kronecker Structure} Now, let us test our approach in the context of the numerical
solution of a 3D reaction--diffusion linear partial differential
equation
\begin{equation}\label{eq:pde_example_forexpA}
\partial_t u = - k \nabla^2 u + \gamma(x,y,z) u.
\end{equation}
Discretizing \eqref{eq:pde_example_forexpA} in the space variables
with second order centered differences, the reference solution can
be computed by means of the matrix $\exp(A)$. We take $k = 1e-8$,
and the action of $\gamma(x,y,z)$ is given by the matrix--vector
product on the semidiscrete equation between $G =
\operatorname{sparsify}(\operatorname{rand}(n^3,n^3))$ where $n$ is
the number of mesh points along one direction of the domain $\Omega
= [0,1]^3$. Function $\operatorname{sparsify}$ gives a sparse
version of $G$ with $0.1\%$ of fill--in and the Laplacian is
discretized with the standard 7-points stencil with homogeneous
Dirichlet conditions, i.e., the semidiscrete equation reads as
\begin{equation*}
\mathbf{u}_t(t) = (A + G)\mathbf{u}(t).
\end{equation*}
The results of this experiment are reported in Table
\ref{tab:exp_pde1_expA}. The reference matrix is computed by using
the incomplete inverse LDU factorization (INVT) that needs two drop
tolerances, $\tau_L = 1e-6$ and $\tau=\tau_Z = 1e-8$. The former is
the drop tolerance for the incomplete $LU$ (or $ILU$ for short)
process and the latter for the post-sparsification of the inversion
of $LU$ factors, respectively; see~\cite{Bertaccini.Filippone12} for
details on approximate inverse preconditioners with inversion of an
$ILU$. A tridiagonal approximation of the correction matrix $\tilde
E = \tilde Z^T \tilde W$ is used, i.e., $\tilde E=g_1(\tilde Z^T
\tilde W)$.
\begin{table}[hbt]
\centering \caption{Execution time in seconds for $\exp(A)$ and
relative errors ($\varepsilon_{\text{rel}}$) with respect to
{\tt{expm(A)}} where $A$ is the discretization matrix of
\eqref{eq:pde_example_forexpA} (the time needed for building the
reference matrix is not considered). INVT with $\tau_L = 1e-6$ and
$\tau=\tau_Z = 1e-8$ is used.}\label{tab:exp_pde1_expA}
\begin{tabular}{ccccccc}
\toprule
& \multicolumn{2}{c}{Direct} & \multicolumn{2}{c}{Update} &
\verb|expm|$(A)$ & \\
$n^3$ & T(s) & $\varepsilon_{\text{rel}}$ & T(s) &
$\varepsilon_{\text{rel}}$ & T(s) & Fill-in \\
\midrule
512 & 0.15 & 2.85e-07 & \textbf{0.07} & 2.82e-07 & 0.92 & 100.00
\% \\
1000 & 0.83 & 2.85e-07 & \textbf{0.35} & 2.83e-07 & 8.19 & 100.00
\% \\
1728 & 4.28 & 2.85e-07 & \textbf{0.94} & 2.83e-07 & 46.23 & 92.40
\% \\
4096 & 118.39 & 2.85e-07 & \textbf{3.72} & 2.84e-07 & 669.39 &
51.77 \% \\
8000 & 834.15 & 2.85e-07 & \textbf{9.69} & 2.82e-07 & 4943.73 &
28.84 \% \\
\bottomrule
\end{tabular}
\end{table}

\subsection{Approximating $\Psi(A)\mathbf{v}$}

\subsubsection{$\exp(A)\mathbf{v}$ - Exponential Decay} To apply our approximation for $\Psi(A)\mathbf{v}$, where $A$ is
large and/or localized and/or possibly structured, we use a Krylov
iterative solver for the systems $(A+\xi_j I)\mathbf{x}=\mathbf{v}$ in \eqref{eq:PFE0}
with and without preconditioning (the corresponding columns will be
labeled as {\emph{Prec}} and {\emph{Not prec}}). The iterative
solvers considered are BiCGSTAB and CG (the latter for symmetric
matrices). The preconditioner is based on the matrix $\tilde Z
(\tilde D+\xi_j \tilde E)^{-1} \tilde W^H$ as in \eqref{eq:Pk}. The
matrix $A$ has the entries as in \eqref{eq:TestMat} while
$\mathbf{v}$ is the normalized unit vector.
\begin{table}[htb]
\centering
\caption{Iterates average  and execution time in seconds for
$\log(A)\mathbf{v}$ for $A$ as in \eqref{eq:TestMat} with
$\alpha=0.2,\,\beta=0.5$. The linear systems are solved with the
Matlab's implementation of {\tt BiCGStab} with and without
preconditioning. INVT with $\tau = \tau_L = \tau_Z = 1e-1$ is used.} \label{tab:itsLogAv}
\begin{tabular}{ccccccc}
\toprule
& \multicolumn{2}{c}{Prec} & \multicolumn{2}{c}{Not prec} \\
n  & iters & T(s) & iters & T(s) \\
\midrule
500&  2 & \textbf{0.05}  & 21 & 0.11 \\
1000& 2 & \textbf{0.05}  & 19 & 0.18  \\
2000&  2 & \textbf{0.08}  & 18 & 0.33  \\
4000& 2 & \textbf{0.95}  & 17 & 2.96  \\
\bottomrule
\end{tabular}
\end{table}
The average of the iterates  in Table \ref{tab:itsLogAv} is much
smaller when the preconditioner is used. Moreover, preconditioned
iterations are independent on the size of the problem. 

In Table \ref{tab:errorExpAv} we report the error, with  respect to
the Matlab's {\tt{expm(A)\textbf{v}}}, of the approximations given by the
{\emph{Prec}} and {\emph{Not prec}} options. The entries in the test
matrix have so a fast decay, since $\alpha=\beta=6$, that the term
$\tilde E$ can be chosen diagonal.
Interestingly, a good accuracy is
reached with respect to the true solution. Moreover, the timings for
the {\emph{Prec}} approach is negligible with respect to that for
the {\emph{Not prec}}.

\begin{table}[htb]
\centering
\caption{Error for $\exp(A)\mathbf{v}$ for $A$ as in \eqref{eq:TestMat} with
$\alpha=\beta=6$. {\emph{Prec}}: our technique with $\tilde E$
diagonal. {\emph{Not prec.}}: \eqref{eq:PFE0} when the linear
systems are solved with the Matlab's {\tt PCG} used without
preconditioning. INVT with $\tau = \tau_L = \tau_Z = 1e-1$ is used.} \label{tab:errorExpAv}
\begin{tabular}{ccc}
\toprule
n  & \textrm{Prec} & \textrm{Not prec} \\
& $\varepsilon_{rel}$ & $\varepsilon_{rel}$ \\
\midrule
500 & 4.5e-06 & 1.8e-08 \\
1000 & 4.5e-06 & 1.8e-08 \\
2000 & 4.5e-06 & 1.8e-08 \\
4000 & 4.5e-06 & 1.8e-08 \\
\bottomrule
\end{tabular}
\end{table}

\subsubsection{$\exp(A)\mathbf{v}$ - Transition Matrices} Let us consider a series of tests matrices of a different nature:
the infinitesimal generators, i.e., transition rate matrices from
the \emph{MARCA} package by Stewart~\cite{william1991marca}. They
are large non--symmetric ill--conditioned matrices whose condition
number ranges from $10^{17}$ to $10^{21}$ and their eigenvalues are
in the square in the complex plane given by $[-90,5.17e-15]\times
i[-3.081,3.081]$. As a first example, we consider the \emph{NCD}
model. It consists of a set of terminals from which the same number
of users issue commands to a system made by a central processing
unit, a secondary memory device and a filling device. In Table
\ref{tab:ncd_expav} we report results for various $n$, obtained by
changing the number of terminals/users. The matrices $A$ are used to
compute $\exp(A)\mathbf{v}$, $\mathbf{v} = (1 \ 2\ \ldots n)^T /n$.
We compare the performance of BiCGSTAB for solving the linear
systems in \eqref{eq:PFE0} without preconditioner and with our
updating strategy, where $g$ extracts only the main diagonal, i.e.,
$g(\cdot)=g_0(\cdot)$. The INVT algorithm with $\tau_Z = 1e-4$ and
$\tau_L = 1e-2$ is used to produce the approximate inverse
factorization. The comparisons consider the time needed for solving
each linear system, i.e., the global time needed to compute
$\exp(A)\mathbf{v}$. Both methods are set to achieve a relative
residual of $10^{-9}$ and the degree of the Chebyshev rational
approximation is $N = 9$. The column $\varepsilon_{\text{rel}}$
reports the relative error between our approximation and
$\operatorname{expm}(A)\mathbf{v}$. For the case with the largest
size $\operatorname{expm}$ gives \lq\lq out of memory\rq\rq error.

\begin{table}[hbt]
\centering
\caption{Approximation of $\exp(A) v$, $A$ from NCD queuing network
example. BiCGSTAB, $N = 9$, tol$= 1e-9$, INVT algorithm with $\tau_Z = 1e-4$ and
$\tau_L = 1e-2$ is used. A $\dagger$ is reported on the $\varepsilon_{\text{rel}}$ when the $\operatorname{expm}$ gives out of memory error and no reference solution is available.}
\label{tab:ncd_expav}
\begin{tabular}{cccccc}
\toprule
& \multicolumn{2}{c}{Not prec} &
\multicolumn{2}{c}{Update} \\
n & iters & T (s) & iters & T (s) & $\varepsilon_{\text{rel}}$ \\
\midrule
286 & 7.50 & 7.62e-03 & 7.50 & \textbf{7.60e-03}  & 8.73e-09 \\
1771 & 17.60 & \textbf{2.95e-02} & 17.60 & 3.00e-02 & 3.46e-07 \\
5456 & 29.00 & \textbf{1.15e-01} & 29.00 & 1.18e-01 & 5.23e-06 \\
8436 & 34.50 & 2.44e-01 & 28.00 & \textbf{1.67e-01} & 1.50e-05 \\
12341 & 43.10 & 3.39e-01 & 33.20 & \textbf{2.68e-01} & 3.87e-05 \\
23426 & 64.30 & 9.66e-01 & 42.30 & \textbf{6.26e-01} & $\dagger$ \\
\bottomrule
\end{tabular}
\end{table}
\subsubsection{$\exp(A)\mathbf{v}$ - Network Adjacency Matrix} 
\begin{figure}[htb]
\centering
\includegraphics[width=0.5\textwidth]{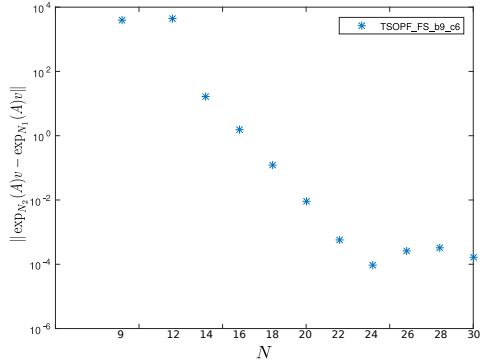}
\caption{Accuracy for various values of $N$ for the
TSOPF\_FS\_b9\_c6 matrix}\label{fig:matrixcollection_expav}
\end{figure}

Let us compute $\exp(A)\mathbf{v}$ for the matrix TSOPF\_FS\_b9\_c6
of dimension $14454$ coming from~\cite{davis2011university}. Results
are reported in Table \ref{tab:matrixcollection_expav} and Figure
\ref{fig:matrixcollection_expav}. For this case we do not have a
reference solution and then the error because MATLAB's \verb|expm|
gives out of memory error. Instead, we consider the Euclidean norm
of the difference of the solutions obtained for consecutive
values of $N$ for $N=6,\ldots,30$. The other settings for the solver
remain unchanged in order to evaluate the efficiency of the
algorithm for the same level of accuracy, i.e., we are using again
the INVT algorithm with $\tau_Z = 1e-4$ and $\tau_L = 1e-2$.

\begin{table}[hbt]
\centering \caption{Approximation of $\exp(A) v$ as the degree $N$
of the Chebyshev approximation varies. The matrix $A$ is
TSOPF\_FS\_b9\_c6~\cite{davis2011university}. The INVT algorithm
with $\tau_Z = 1e-4$ and $\tau_L = 1e-2$ is used.}
\label{tab:matrixcollection_expav} {
\begin{tabular}{ccccc}
\toprule
\multicolumn{5}{c}{Matrix TSOPF\_FS\_b9\_c6}\\
\multicolumn{5}{c}{Size: 14454, $\kappa_2(A)=$3.1029e+12} \\
\midrule
\multicolumn{2}{c}{Not prec} &
\multicolumn{2}{c}{Prec} \\
iters & T(s) & iters & T(s) & $N$ \\
\midrule
171.33 & 1.22e+00 & 15.00 & \textbf{2.62e-01} & 6 \\
145.20 & 1.38e+00 & 36.50 & \textbf{1.05e+00} & 9 \\
99.58 & 1.44e+00 & 8.75 & \textbf{2.33e-01} & 12 \\
77.93 & 1.32e+00 & 7.64 & \textbf{2.29e-01} & 14 \\
70.38 & 1.25e+00 & 7.06 & \textbf{2.44e-01} & 16 \\
71.00 & 1.45e+00 & 6.61 & \textbf{2.65e-01} & 18 \\
59.70 & 1.34e+00 & 6.15 & \textbf{2.80e-01} & 20 \\
53.32 & 1.32e+00 & 5.95 & \textbf{2.93e-01} & 22 \\
51.67 & 1.38e+00 & 5.75 & \textbf{3.14e-01} & 24 \\
46.65 & 1.37e+00 & 5.58 & \textbf{3.30e-01} & 26 \\
44.86 & 1.44e+00 & 5.39 & \textbf{3.49e-01} & 28 \\
43.30 & 1.47e+00 & 5.20 & \textbf{3.66e-01} & 30 \\
\bottomrule
\end{tabular}}
\end{table}
We observe two different effects for higher degree of
approximations in Table~\ref{tab:matrixcollection_expav}. On one hand, from Figure~\ref{fig:matrixcollection_expav},
the relative error is reduced, as
expected from the theoretical analysis, while, on the other, it
makes the shifted linear system more well--conditioned. Note that the
gain obtained using our preconditioning strategy is sensible even
for large matrices.

\subsubsection{$\log(A)\mathbf{v}$ - Polynomial Decay} Let us consider the matrix $A$ from~\cite{lu1998computing} with
entries given by
\begin{equation}\label{eq:log_exp1}
a_{i,j} = \frac{1}{2+(i-j)^2},\quad i,j=1,\ldots,n
\end{equation}
in order to approximate $\log(A) \mathbf{v}$, $\mathbf{v} = (1,1,\ldots,1)^T$. $A$
is symmetric positive definite with a minimum eigenvalue of the
order of $10^{-2}$ and its entries decay polynomially. We
approximate $\log(A)\mathbf{v}$ with \eqref{eq:PFE0}; BiCGSTAB is
used with our preconditioner update strategy and without it (Not
prec). The seed preconditioner is computed using INVT with  $\tau_Z
= 1e-1$ and $\tau_L = 1e-2$. We include the results with MATLAB's
$\operatorname{logm}(A)\mathbf{v}$. In particular, we use $N = 30$
for the approximation of the logarithm function. Results are
collected in Table~\ref{tab:logav_exp1decay}.
\begin{table*}[hbt]
\centering
\caption{Computation of $\log(A)\mathbf{v}$ with $A$ as in equation
\eqref{eq:log_exp1}. Note the moderate decay and a spectrum
that ranges in the interval $[6e-2,3]$. For INVT $\tau_Z = 1e-1$ and $\tau_L = 1e-2$ are used.}\label{tab:logav_exp1decay}
\begin{tabular}{cccccccc}
\toprule
BiCGSTAB & \multicolumn{2}{c}{Not prec} & \multicolumn{3}{c}{Update} & \multicolumn{2}{c}{$\operatorname{logm}(A)\mathbf{v}$}\\
n & iters & T(s) &  iters & T(s)
& Fill--In & T (s) & $\varepsilon_{\text{rel}}$ \\
\midrule
1000 & 11.88 & 2.4 & 5.07 &
1.25 & 3.16 \% & \textbf{0.169} & 1.91e-06 \\
4000 & 11.05 & 34.2 & 4.73 &
16.8 & 0.80 \% & \textbf{14.7} & 1.54e-06 \\
8000 & 10.58 & 133.78 & 4.53 &
\textbf{65.7} & 0.40 \% & 116.9 & 1.66e-06 \\
12000 & 10.32 & 296.9 & 4.38 &
\textbf{145.7} & 0.27 \% & 428.2 & 1.74e-06 \\
\bottomrule
\end{tabular}
\end{table*}

\subsubsection{$\log(A)\mathbf{v}$ - Matrix Collection} Finally, we consider some matrices from \emph{The University of
Florida Sparse Matrix Collection} (see~\cite{davis2011university}),
focusing on INVT with a seed  preconditioner with $\tau_Z = 1e-1$
and $\tau_L = 1e-2$.
The results are collected in Table \ref{tab:logav_florida}, 
and confirm what we observed in the other tests.
\begin{table*}[h]
\centering \caption{Approximation of $\log(A)\mathbf{v}$ with $A$
SPD from The University of Florida Sparse Matrix Collection. The
real parts of the eigenvalues are all in the interval
$[2.324e-14,1.273e+08]$. The updated preconditioners are computed
using INVT with $\tau_Z = 1e-1$ and $\tau_L
=1e-2$.}\label{tab:logav_florida} {
\begin{tabular}{ccccccccc}
    \toprule
    \multicolumn{2}{c}{BiCGSTAB} &  \multicolumn{2}{c}{Not prec} & \multicolumn{3}{c}{Update} & \multicolumn{2}{c}{$\operatorname{logm}(A)\mathbf{v}$}\\
    Name & n & iters & T(s) &  iters
    & T(s) & Fill--In & T(s) & $\varepsilon_{\text{rel}}$ \\
    \midrule
    1138\_bus & 1138 & 198.93 & 1.3 & 31.18 & 0.4 & 0.84 \%
    &
    \textbf{0.22} & 4.41e-07 \\
    Chem97ZtZ & 2541 & 27.98 & 0.34 & 6.43 & \textbf{0.12} &
    0.10 \% & 3.5 & 1.87e-07 \\
    bcsstk21 & 3600 & 157.85 & 4.8 & 76.4 & \textbf{3.1} &
    1.36 \% & 10 & 3.10e-07 \\
    t2dal\_e & 4257 & 232.00 & 4.2 & 98.90 & \textbf{1.78} &
    0.02 \% & 2.58 & 6.82e-04 \\
    crystm01 & 4875 & 23.35 & 1.03 & 11.48 & \textbf{0.56} &
    0.17 \% & 25.3 & 3.16e-07 \\
    \bottomrule
\end{tabular}}

\end{table*}

\subsection{$\Psi(A)\mathbf{v}$ with updates and with Krylov subspace methods}

A popular class of effective algorithms for approximating
$\Psi(A)\mathbf{v}$ for a given large and sparse matrix $A$ relies
on Krylov subspace methods. The basic idea is to project the problem
into a smaller space and then to make its solution potentially
cheaper. The favorable computational and approximation properties
have made the Krylov subspace methods extensively used; see, e.g.,
\cite{Saad92,LopezSimoncini06,Hochbruck.Lubich.97,Moret2009}.

Over the  years, some tricks have been added to these techniques to
make them more effective, both in terms of computational cost and
memory requirements, see, e.g.,
\cite{Moret.Novati.04b,vandenEshof.Hochbruck.06,Popolizio.Simoncini.08,AfanasjewEtAl08,KnizhnermanSimoncini2010,MoretPopolizio12}.
In particular, as shown by Hochbruck and Lubich
\cite{Hochbruck.Lubich.97}, the convergence depends on the spectrum
of $A$. For our test matrices the spectrum has just a moderate
extension in the complex plane. Thus, the underlying Krylov subspace
techniques for approximating $\Psi(A)\mathbf{v}$ can be appropriate.

The approximation spaces for these techniques are defined as
$$
{\cal{K}}_m(A,\mathbf{v})={\textrm {span}}\{\mathbf{v},A\mathbf{v},\ldots,A^{m-1}\mathbf{v}\}.
$$
Since the basis given by the vectors
$\mathbf{v}=\mathbf{v}_1,A\mathbf{v}=\mathbf{v}_2,\ldots,A^{m-1}\mathbf{v}=\mathbf{v}_m$ can be very ill--conditioned, one
usually applies the modified Gram-Schmidt method to get an
orthonormal basis with starting vector $\mathbf{v}_1=\mathbf{v}/\|\mathbf{v}\|$. Thus, if these
vectors $\mathbf{v}_1,\ldots,\mathbf{v}_m$ are the columns of a matrix $V_m$ and the
upper Hessenberg matrix $H_m$ collects the coefficients $h_{i,j}$ of
the orthonormalization process, the following expression by Arnoldi
holds
$$
AV_m = V_m H_m+h_{m+1,m} \mathbf{v}_{m+1}\mathbf{e}_m^T,
$$
where $\mathbf{e}_m$ denotes the $m$th column of the identity matrix. An
approximation to $\Psi(A)\mathbf{v}$ can be obtained as
$$
y_m=\|\mathbf{v}\| V_m \Psi(H_m)\mathbf{e}_1.
$$
The procedure reduces to the three-term Lanczos  recurrence when $A$
is symmetric, which results in a tridiagonal matrix $H_m$. One has
still to face the issue of evaluating a matrix function, but, if
$m\ll n$, for the matrix $H_m$, which is just $m\times m$. Several
approaches can then be tried. For example, one can use the built--in
function $\tt{funm}$ in Matlab, based on the Schur decomposition of
the matrix argument, and the Schur-Parlett algorithm to evaluate the
function of the triangular factor \cite{Higham2008}.

We consider the application of our strategy for the computation of
$\exp(A)\mathbf{v}$, with a matrix $A$ generated from the discretization
of the following 2D advection-diffusion problem 
\begin{equation}\label{eq:advectiondiffusionproblem}
\left\lbrace\begin{array}{rl}
u_t = & \displaystyle\frac{\partial }{\partial x}\left(k_1\frac{\partial
u}{\partial x}\right) + \frac{\partial }{\partial
y}\left(k_2(y)\frac{\partial
u}{\partial y}\right) \\ &\displaystyle + t_1(x) \frac{\partial u}{\partial x}  + t_2(y)
\frac{\partial u}{\partial y}, \quad x \in [0,1]^2   \\
u(x,y,t) = & 0, \quad x \in \partial [0,1]^2, \\
u(x,y,0) = & u_0(x,y)
\end{array}\right.
\end{equation}
where the coefficients are $k_1 = 10^{-2}$, $k_2(x) = 2+10^{-5}
\cos(5\pi x)$, $t_1(x) = 1+0.15\sin(10\pi x)$ and $t_2(x) =
1+0.45\sin(20\pi x)$. The second order centered differences and
first order upwind are used to discretize the Laplacian and the
convection terms, respectively. The purpose of this experiment,
whose results are in Table \ref{tab:lanczos}, is comparing the
performance of our updating approach, using INVT with $\tau_L =
1e-5$ and $\tau_Z = 1e-2$, with a Krylov subspace method. For the
latter we use the classical stopping criterion based on monitoring
$$
\gamma=h_{m+1,m}|\mathbf{e}_m^T \exp(H_m)\mathbf{e}_1|.
$$
We stop the iteration when $\gamma$ becomes smaller than
$10^{-6}$. The threshold $\gamma$ was tuned to the accuracy expected
by the {\emph{Update}} approach.

\begin{table}[htb]
\centering
\caption{Errors and execution time for $\exp(A)\mathbf{v}$ for
    $A$ obtained as the finite difference discretization of
    \eqref{eq:advectiondiffusionproblem}. INVT with $\tau_L = 1e-5$ and
    $\tau_Z = 1e-2$ is used.} \label{tab:lanczos}
{   
    \begin{tabular}{ccccc}
        \toprule
        n  & \multicolumn{2}{c}{\textrm{Update}} & \multicolumn{2}{c}{\textrm{Arnoldi}}  \\
        & $\varepsilon_{\text{rel}}$ & T(s) & $\varepsilon_{\text{rel}}$ & T(s) \\
        \midrule
        100 & 2.64e-06 & 6.74e-02  & 3.51e-06 & \textbf{2.91e-02} \\
        196 & 3.81e-06 & \textbf{7.80e-02} & 1.20e-06 & 1.09e-01 \\
        484 & 1.22e-08 & \textbf{2.97e-01} & 1.60e-08 & 5.23e-01 \\
        961 & 7.68e-07 & \textbf{1.27e+00} & 1.68e-07 & 2.70e+00 \\
        \bottomrule
    \end{tabular}}
\end{table}
From these experiences and other non reported here, we can conclude
that our techniques, under appropriate hypotheses of sparsity or
locality of the matrices, seem to have reasonably comparable
performances with Krylov methods when computing $\Psi(A)\mathbf{v}$.

Moreover, we can expect even more interesting performances when
simultaneous computations of vectors such as
$\Psi(A)\mathbf{w}_1,\Psi(A)\mathbf{w}_2,\ldots,\Psi(A)\mathbf{w}_K$
are required. This will give us another level of parallelism beyond
the one that can be exploited in the simultaneous computation of the
terms in \eqref{eq:PFE0}. In particular, this can be true when $K$
is large and each vector $\mathbf{w}_{j}$ does depend on the
previous values $\mathbf{w}_i$ and $\Psi(A)\mathbf{w}_i$. In this
setting we can construct the factors once in order to reduce the
impact of the initial cost for computing the approximate inverse
factors. A building cost that can be greatly reduced by using
appropriate parallel algorithms and architectures;
see~\cite{Bertaccini.Filippone12} and references therein.

\section{Conclusions}
\label{sec:Discussion}

Consider the hypotheses used in this research:
\begin{itemize}
    \item the function $\Psi$ should be smooth enough in the sense of \cite{HaleHighamTref08};
    \item the off--diagonal entries of $A^{-1}$ should decay fast enough away from the main diagonal.
\end{itemize}
A natural question arises: which is the most relevant feature to
make our {\emph{Update}} approach effective? Several papers have
been devoted to the analysis of the decay of the entries of
$\Psi(A)$ when  the behavior of the entries of $A$ is known
\cite{Demko.Moss.Smith84,BenziGolub99,Iserles00,Benzi.Razouk.07}. A
unifying analysis has been proposed in~\cite{Benzi.Razouk.07}, where
the influence of $A$ and $\Psi$ is considered. One of their results
shows that the entries in $\Psi(A)$ can be bounded by a constant and
a term depending on the decay of the entries of $A^{-1}$, provided
that $\Psi$ is analytic in a suitable region containing the spectrum
of~$A$.

As an example, let us recall a result concerning band matrices
(see \cite[Corollary 3.6]{Benzi.Razouk.07}). If
$A$ is a diagonalizable band matrix of dimension $n\times n$, then
we have
$$
|\Psi(A)_{i,j}|< c \, \kappa_2(X) \lambda^{|i-j|},\,\,i,j=1\ldots,n
$$
where $X$ is the matrix of eigenvectors of $A$, $\kappa(\cdot)$ is
the condition number in the Euclidean norm, $c$ is a positive
constant and $0<\lambda<1$. We suppose that $\kappa_2(X)$ is
moderate, otherwise the above bound is useless.

The latter bound seems to give a slightly greater importance to the
matrix argument $A$ than to the function $\Psi$.

This discussion confirms our experience: the performances of the
{\emph{Update}} approach seem to depend more on $A$ (in particular,
on the behavior of the entries of $A^{-1}$ away from the main
diagonal) than on $\Psi$. Indeed, as can be observed in
Table~\ref{tab:errorExpMat}, if we approximate $\exp(A)$, for three
different matrix arguments with different decays, the accuracy
changes seem to be more influenced by the matrix argument. A similar
comment can be done also for the plots in Figure~\ref{ErrorVsBand}.

\section*{Funding} This work was supported in part by
INDAM-GNCS 2018 projects ``Tecniche innovative per problemi di
algebra lineare'' and ``Risoluzione numerica di equazioni di evoluzione
integrali e differenziali con memoria''.

\bibliographystyle{tfs}
\bibliography{BibShiftSyst}

\end{document}